\documentclass[11pt]{amsart}
 \usepackage{amssymb,amsfonts,amscd,amsbsy}
\usepackage[mathscr]{eucal}
\usepackage{url}

\theoremstyle{plain}
\newtheorem{thm}{Theorem}[section]
\newtheorem{lem}[thm]{Lemma}
\newtheorem{prop}[thm]{Proposition}

\theoremstyle{definition}

\numberwithin{equation}{section}

\begin{document}
\title[Rank differences for overpartitions]{Rank differences for overpartitions}
\author{Jeremy Lovejoy and Robert Osburn}

\address{CNRS, LIAFA, Universit{\'e} Denis Diderot, 2, Place Jussieu, Case 7014, F-75251 Paris Cedex 05, FRANCE}

\address{School of Mathematical Sciences, University College Dublin, Belfield, Dublin 4, Ireland}

\address{IH{\'E}S, Le Bois-Marie, 35, route de Chartres, F-91440 Bures-sur-Yvette, FRANCE}

\email{lovejoy@liafa.jussieu.fr}

\email{robert.osburn@ucd.ie, osburn@ihes.fr}

\thanks{The first author was partially supported by an ACI ``Jeunes Chercheurs et Jeunes Chercheuses".}


\date{May 28, 2007}

\begin{abstract}
In 1954, Atkin and Swinnerton-Dyer proved Dyson's conjectures on
the rank of a partition by establishing formulas for the
generating functions for rank differences in arithmetic
progressions. In this paper, we prove formulas for the generating
functions for rank differences for overpartitions.  These are in
terms of modular functions and generalized Lambert series.
\end{abstract}

\maketitle
\section{Introduction}
The rank of a partition is the largest part minus the number of
parts.  This statistic was introduced by Dyson \cite{dyson}, who
observed empirically that it provided a combinatorial explanation
for Ramanujan's congruences $p(5n+4) \equiv 0 \pmod{5}$ and
$p(7n+5) \equiv 0 \pmod{7}$.  Here $p(n)$ denotes the usual
partition function.  Specifically, Dyson conjectured that if
$N(s,m,n)$ denotes the number of partitions of $n$ whose rank is
congruent to $s$ modulo $m$, then for all $0 \leq s \leq 4$ and $0
\leq t \leq 6$ we have
$$
N(s,5,5n+4) = \frac{p(5n+4)}{5}
$$
and
$$
N(t,7,7n+5) = \frac{p(7n+5)}{7}.
$$

Atkin and Swinnerton-Dyer proved these assertions in 1954
\cite{asd}. In fact, they proved much more, establishing
generating functions for every rank difference $N(s,\ell,\ell n +
d) - N(t,\ell,\ell n + d)$ with $\ell = 5$ or $7$ and $0 \leq
d,s,t < \ell$.  Many of these turned out to be non-trivially $0$,
while others were infinite products and still others were
generalized Lambert series related to Ramanujan's third order mock
theta functions.  Such formulas in the case $\ell = 11$ were
subsequently given by Atkin and Hussain \cite{At-Hu1} in a similar
(though technically far more difficult) manner.

Dyson's rank extends in the obvious way to overpartitions. Recall
that an overpartition \cite{Co-Lo1} is simply a partition in which
the first occurrence of each distinct number may be overlined. For
example, the $14$ overpartitions of $4$ are

\begin{center}
$4$, $\overline{4}$, $3+1$, $\overline{3} + 1$,
$3 + \overline{1}$, $\overline{3} + \overline{1}$,
$2+2$, $\overline{2} + 2$, $2+1+1$, $\overline{2} + 1 + 1$,
$2+ \overline{1} + 1$, $\overline{2} + \overline{1} + 1$, $1+1+1+1$, $\overline{1} + 1 + 1 +1$.
\end{center}

\noindent Overpartitions (often under other names) naturally arise
in diverse areas of mathematics where partitions already occur,
such as mathematical physics \cite{fjm,fjm1}, symmetric
functions \cite{Bre1,De-La-Ma3}, representation theory
\cite{Ka-Kw1}, and algebraic number theory \cite{love,love2}.
With basic hypergeometric series, Dyson's rank for overpartitions
and its generalizations play an important role in combinatorial
studies of Rogers-Ramanujan type identities
\cite{Co-Lo-Ma1,Co-Ma1,Lo1}.

As for arithmetic properties, let $\overline{N}(s,m,n)$ be the
number of overpartitions of $n$ with Dyson's rank congruent to $s$
modulo $m$. It has recently been shown \cite{Br-Lo1} that if $u$
is any natural number, $m$ is odd and $\ell \geq 5$ is a prime
such that $\ell \nmid 6m$ or $\ell^j = m$, then there are
infinitely many non-nested arithmetic progressions $An+B$ such
that
$$
\overline{N}(s,m,An+B) \equiv 0 \pmod{\ell^u}
$$
for all $0 \leq s < m$.  This is analogous to congruences
involving Dyson's rank for partitions \cite{Br-On1}.

On the other hand, there are no congruences of the form
$\overline{p}(\ell n + d) \equiv 0 \pmod{\ell}$ for primes $\ell
\geq 3$ \cite{Choi1}.  Here $\overline{p}(n)$ denotes the number
of overpartitions of $n$. The generating functions for the rank
differences $\overline{N}(s,\ell,\ell n+d) -
\overline{N}(t,\ell,\ell n + d)$ then provide a measure of the
extent to which the rank fails to produce a congruence
$\overline{p}(\ell n + d) \equiv 0 \pmod{\ell}$.  In this paper,
we find formulas for these generating functions for $\ell = 3$ and
$5$ in terms of modular functions and generalized Lambert series.
Using the notation
\begin{equation} \label{rst}
R_{st}(d) = \sum_{n \geq 0} \left(\overline{N}(s,\ell,\ell n + d)
- \overline{N}(t,\ell,\ell n + d)\right ) q^n,
\end{equation}
where the prime $\ell$ will always be clear, the main results are
summarized in Theorems \ref{main3} and \ref{main} below.

\begin{thm}\label{main3}
For $\ell=3$, we have

\begin{equation} \label{r_01(0)}
R_{01}(0) = -1 +
\frac{(q^3;q^3)_{\infty}^2(-q;q)_{\infty}}{(q)_{\infty}(-q^3;q^3)_{\infty}^2},
\end{equation}

\begin{equation} \label{r_01(1)}
R_{01}(1) =
\frac{2(q^3;q^3)_{\infty}(q^{6};q^{6})_{\infty}}{(q)_{\infty}},
\end{equation}

\begin{equation}
R_{01}(2) =
\frac{4(-q^3;q^3)_{\infty}^2(q^{6};q^{6})_{\infty}^2}{(q^2;q^2)_{\infty}}
- \frac{6(-q^3;q^3)_{\infty}}{(q^3;q^3)_{\infty}} \sum_{n \in
\mathbb{Z}} \frac{(-1)^nq^{3n^2+3n}}{1-q^{3n+1}}.
\end{equation}
\end{thm}

\begin{thm} \label{main} For $\ell=5$, we  have

\begin{equation} \label{r_{12}(0)}
R_{12}(0)=\frac{2q(q^{10}; q^{10})_{\infty}}{(q^{3}, q^{4}, q^{6},
q^{7}; q^{10})_{\infty}},
\end{equation}

\begin{equation} \label{r_{12}(1)}
R_{12}(1)= \frac{-2q(-q^{5}; q^{5})_{\infty}}{(q^{5};
q^{5})_{\infty}} \sum_{n \in \mathbb{Z}}
\frac{(-1)^nq^{5n^2+5n}}{1-q^{5n+2}} ,
\end{equation}

\begin{equation} \label{r_{12}(2)}
R_{12}(2)=\frac{2(q^{10}; q^{10})_{\infty}}{(q, q^{4};
q^{5})_{\infty}},
\end{equation}

\begin{equation} \label{r_{12}(3)}
R_{12}(3)=\frac{-2(q^{10}; q^{10})_{\infty}}{(q^{2}, q^{3};
q^{5})_{\infty}},
\end{equation}

\begin{equation} \label{r_{12}(4)}
R_{12}(4) = \frac{6(-q^{5}; q^{5})_{\infty}}{(q^{5};
q^{5})_{\infty}}\sum_{n \in
\mathbb{Z}}\frac{(-1)^nq^{5n^2+5n}}{1-q^{5n+1}} - \frac{4(q^{2},
q^{8},q^{10}; q^{10})_{\infty}} {(q^{4}, q^{6};
q^{10})_{\infty}^{2} (q, q^{9}; q^{10})_{\infty}},
\end{equation}

\begin{equation} \label{r_{02}(0)}
R_{02}(0) = -1  + \frac{(-q^2, -q^{3}; q^{5})_{\infty} (q^{5};
q^{5})_{\infty}}{(q^2, q^{3}; q^{5})_{\infty} (-q^{5};
q^{5})_{\infty}},
\end{equation}

\begin{equation} \label{r_{02}(1)}
R_{02}(1) = \frac{2(q^{4}, q^{6}, q^{10};
q^{10})_{\infty}}{(q^{2}, q^{8};
q^{10})_{\infty}^2(q^3,q^7;q^{10})_{\infty}} + \frac{4q(-q^{5};
q^{5})_{\infty}}{(q^{5}; q^{5})_{\infty}}\sum_{n \in
\mathbb{Z}}\frac{(-1)^nq^{5n^2+5n}}{1-q^{5n+2}}  ,
\end{equation}

\begin{equation} \label{r_{02}(2)}
R_{02}(2)=0,
\end{equation}

\begin{equation} \label{r_{02}(3)}
R_{02}(3)=\frac{2(q^{10}; q^{10})_{\infty}}{(q^{2}, q^{3};
q^{5})_{\infty}},
\end{equation}

\begin{equation} \label{r_{02}(4)}
R_{02}(4)=\frac{2(q^{2}, q^{8}, q^{10}; q^{10})_{\infty}} {(q^{4},
q^{6}; q^{10})_{\infty}^{2} (q, q^{9}; q^{10})_{\infty}} -
\frac{2(-q^{5}; q^{5})_{\infty}}{(q^{5}; q^{5})_{\infty}} \sum_{n
\in \mathbb{Z}}\frac{(-1)^nq^{5n^2+5n}}{1-q^{5n+1}}.
\end{equation}

\end{thm}
Here we have employed the standard basic hypergeometric series
notation \cite{Ga-Ra1},
$$
(a_1,a_2,\dots,a_j;q)_n = \prod_{k=0}^{n-1}(1-a_1q^k)(1-a_2q^k)
\cdots (1-a_jq^k),
$$
following the custom of dropping the ``$;q$" unless the base is
something other than $q$.  We should also remark that if the
number of overpartitions of $n$ with rank $m$ is denoted by
$\overline{N}(m,n)$, then conjugating Ferrers diagrams shows that
$\overline{N}(m,n) = \overline{N}(-m,n)$ \cite{Lo1}. Hence the
values of $s$ and $t$ considered in Theorems \ref{main3} and
\ref{main} are sufficient to find any rank difference generating
function $R_{st}(d)$.

To prove our main theorems we adapt the method of Atkin and
Swinnerton-Dyer \cite{asd}.  This may be generally described as
regarding groups of identities as equalities between polynomials
of degree $\ell-1$ in $q$ whose coefficients are power series in
$q^{\ell}$. Specifically, we first consider the expression

\begin{equation} \label{term}
\sum_{n=0}^{\infty} \Bigl\{ \overline{N}(s,\ell,n) -
\overline{N}(t,\ell,n) \Bigr\} q^{n}
\frac{(q)_{\infty}}{2(-q)_{\infty}}.
\end{equation}

\noindent By (\ref{gen1}), (\ref{s}), and (\ref{final}), we write
(\ref{term}) as a polynomial in $q$ whose coefficients are power
series in $q^{\ell}$.  We then alternatively express (\ref{term})
in the same manner using Theorem \ref{main} and Lemma \ref{lem6}.
Finally, we use various $q$-series identities to show that these
two resulting polynomials are the same for each pair of values of
$s$ and $t$.  We should stress that this technique requires
knowing all of the generating function formulas for the rank
differences beforehand.  We cannot prove a generating function for
$R_{st}(d)$ for some $d$ without proving them for all $d$.

The paper is organized as follows.  In Section 2 we collect some
basic definitions, notations and generating functions.  In Section
3 we record a number of equalities between an infinite product and
a sum of infinite products.  These are ultimately required for the
simplification of identities that end up being more complex than
we would like, principally because there is only one $0$ in
Theorem \ref{main}. In Section 4 we prove two key $q$-series
identities relating generalized Lambert series to infinite
products, and in Section 5 we give the proofs of Theorems
\ref{main3} and \ref{main}.

Before proceeding, it is worth mentioning that there is a
theoretical reason why some of the rank differences in the main
theorems are modular and others are not.  We do not go into great
detail here, but the weak Maass forms that lie behind the rank
differences $\overline{N}(s,\ell,\ell n + d) -
\overline{N}(t,\ell,\ell n + d)$ \cite{Br-Lo1} can be shown in
many cases (such as when $-d$ is not a square modulo $\ell$)
to be weakly holomorphic modular forms.  This
has been carried out in detail for the rank differences for
ordinary partitions in \cite{Br-On-Rh1}.  For now, however, it
seems that the groups involved are too small (and the number of
inequivalent cusps too large) to justify pursuing proofs of
identities for rank differences using this framework. In any case,
the present technique will always have the advantage of providing
formulas for all of the generating functions for rank differences,
modular or not.

\section{Preliminaries}
We begin by introducing some notation and definitions, essentially
following \cite{asd}.  With $y=q^{\ell}$, let

$$
r_{s}(d):=\sum_{n=0}^{\infty} \overline{N}(s,\ell, \ell n+d) y^n
$$
and
$$
r_{st}(d):=r_{s}(d) - r_{t}(d).
$$
Thus we have
$$
\sum_{n=0}^{\infty} \overline{N}(s,\ell,n) q^n =
\sum_{d=0}^{\ell-1} r_{s}(d) q^d.
$$

To abbreviate the sums occurring in Theorems \ref{main3} and \ref{main}, we define
$$
\Sigma(z,\zeta, q):= \sum_{n \in \mathbb{Z}} \frac{(-1)^n
\zeta^{n} q^{n^2 + n}}{1-zq^n}.
$$
Henceforth we assume that $a$ is not a multiple of $\ell$.  We
write
$$
\Sigma(a,b):=\Sigma(y^a, y^b, y^{\ell})=\sum_{n \in \mathbb{Z}}
\frac{(-1)^n y^{bn + \ell n(n+1)}}{1-y^{\ell n+a}}.
$$
and
$$
\Sigma(0,b) :=  \sideset{}{'} \sum_{n \in \mathbb{Z}} \frac{(-1)^n
y^{bn + \ell n(n+1)}}{1-y^{\ell n}},
$$
where the prime means that the term corresponding to $n=0$ is
omitted.

To abbreviate the products occurring in Theorems \ref{main3} and \ref{main}, we
define

$$
P(z,q):=\prod_{r=1}^{\infty} (1-zq^{r-1})(1-z^{-1}q^r),
$$

$$
P(a):=P(y^a, y^{\ell}),
$$

\noindent and

$$
P(0):=\prod_{r=1}^{\infty} (1-y^{\ell r}).
$$

\noindent Note that $P(0)$ is not $P(a)$ evaluated at $a=0$. We also have the relations

\begin{equation} \label{p1}
P(z^{-1}q, q)=P(z,q)
\end{equation}

\noindent and

\begin{equation} \label{p2}
P(zq, q)=-z^{-1}P(z,q).
\end{equation}

\noindent From (\ref{p1}) and (\ref{p2}), we have

\begin{equation} \label{p3}
P(\ell-a)=P(a)
\end{equation}

\noindent and

\begin{equation} \label{p4}
P(-a)=P(\ell+a)=-y^{-a} P(a).
\end{equation}

In \cite{Lo1}, it is shown that the two-variable generating
function for $\overline{N}(m,n)$ is

\begin{equation} \label{gen}
\sum_{n=0}^{\infty} \overline{N}(m,n) q^n = \frac{2(-q)_{\infty}}
{(q)_{\infty}} \sum_{n=1}^{\infty} (-1)^{n-1} q^{n^2 + |m|n}
\frac{1-q^n}{1+q^n}.
\end{equation}
From this we may easily deduce that the generating function for
$\overline{N}(s,m,n)$ is

\begin{equation} \label{gen1}
\sum_{n=0}^{\infty} \overline{N}(s,m,n) q^n =
\frac{2(-q)_{\infty}}{(q)_{\infty}} \sideset{}{'} \sum_{n \in
\mathbb{Z}} \frac{(-1)^nq^{n^2 + n}(q^{sn} +
q^{(m-s)n})}{(1+q^n)(1 - q^{mn})}.
\end{equation}
Hence it will be beneficial to consider sums of the form

\begin{equation} \label{s}
\overline{S}(b):= \sideset{}{'}\sum_{n \in \mathbb{Z}}
\frac{(-1)^nq^{n^2+bn}}{1-q^{\ell n}}.
\end{equation}
We will require the relation

\begin{equation} \label{rels}
\overline{S}(b)=-\overline{S}(\ell-b),
\end{equation}
which follows from the substitution $n \to -n$ in (\ref{s}). We
shall also exploit the fact that the functions
$\overline{S}(\ell)$ are essentially infinite products.

\begin{lem}\label{Sofq}
We have
$$
\overline{S}(\ell) = \frac{-(q)_{\infty}}{2(-q)_{\infty}} +
\frac{1}{2}.
$$
\end{lem}

\begin{proof}
Using the relation \eqref{rels}, we have
\begin{eqnarray*}
-2 S(\ell) &=& -2\sideset{}{'} \sum_{n \in \mathbb{Z}}
\frac{(-1)^nq^{n^2 + \ell n}}{1-q^{\ell n}} \\ &=& \sideset{}{'}
\sum_{n \in \mathbb{Z}} \frac{(-1)^nq^{n^2}}{1-q^{\ell n}} -
\sideset{}{'} \sum_{n \in \mathbb{Z}} \frac{(-1)^nq^{n^2 + \ell
n}}{1-q^{\ell n}} \\
&=& \sideset{}{'} \sum_{n \in \mathbb{Z}} (-1)^nq^{n^2}.
\end{eqnarray*}
The lemma now follows upon applying the case $z = -1$ of Jacobi's
triple product identity,
\begin{equation} \label{jtp}
\sum_{n \in \mathbb{Z}} z^nq^{n^2} = (-zq,-q/z,q^2;q^2)_{\infty}.
\end{equation}
\end{proof}

\section{Infinite product identities}
In this section we record some identities involving infinite
products.  These will be needed later on for simplification and
verification of certain identities. First, we have a result which
is the analogue of Lemma 6 in \cite{asd}.

\begin{lem} \label{lem6} We have

\begin{equation} \label{lem6eq1}
\frac{(q)_{\infty}}{(-q)_{\infty}} = \frac{(q^{9};
q^{9})_{\infty}}{(-q^{9}; q^{9})_{\infty}}
 - 2q(q^{3}, q^{15}, q^{18}; q^{18})_{\infty}
\end{equation}
and
\begin{equation} \label{lem6eq2}
\frac{(q)_{\infty}}{(-q)_{\infty}} = \frac{(q^{25};
q^{25})_{\infty}}{(-q^{25}; q^{25})_{\infty}}
 - 2q(q^{15}, q^{35}, q^{50}; q^{50})_{\infty} + 2q^4 (q^{5}, q^{45}, q^{50}; q^{50})_{\infty}.
\end{equation}

\end{lem}

\begin{proof}
This really just amounts to two special cases of \cite[Theorem
1.2]{An-Hi1}.  Indeed, \eqref{lem6eq1} is \cite[Eq.
(1.18)]{An-Hi1}. We give the details for \eqref{lem6eq2}.
Beginning with Jacobi's triple product identity \eqref{jtp}, we
have
\begin{eqnarray*}
\frac{(q)_{\infty}}{(-q)_{\infty}} &=& \sum_{n \in \mathbb{Z}
\atop n \equiv 0 \pmod{5}} (-1)^nq^{n^2} + \sum_{n \in \mathbb{Z}
\atop n \equiv \pm 1 \pmod{5}} (-1)^nq^{n^2} + \sum_{n \in
\mathbb{Z} \atop n \equiv \pm 2 \pmod{5}} (-1)^nq^{n^2} \\
&=& \sum_{n \in \mathbb{Z}}(-1)^nq^{25n^2} + 2\sum_{n \in
\mathbb{Z}}(-1)^{5n+1}q^{(5n+1)^2} + 2\sum_{n \in
\mathbb{Z}}(-1)^{5n+2}q^{(5n+2)^2} \\
&=& \sum_{n \in \mathbb{Z}}(-1)^nq^{25n^2} - 2q\sum_{n \in
\mathbb{Z}}(-1)^{n}q^{25n^2+10n} + 2q^4\sum_{n \in
\mathbb{Z}}(-1)^{n}q^{25n^2+20n}.
\end{eqnarray*}
Again using \eqref{jtp}, we obtain the right hand side of
\eqref{lem6eq2}.
\end{proof}

Next, we quote a result of Hickerson \cite[Theorem 1.1]{Hi1} along
with some of its corollaries.
\begin{lem} \label{Hick1}
$$
P(x,q)P(z,q)(q)_{\infty}^2 =
P(-xz,q^2)P(-qz/x,q^2)(q^2;q^2)_{\infty}^2 -
xP(-xzq,q^2)P(-z/x,q^2)(q^2;q^2)_{\infty}^2.
$$
\end{lem}

The first corollary was recorded by Hickerson \cite[Theorem
1.2]{Hi1}.  It follows by applying Lemma \ref{Hick1} twice, once
with $x$ replaced by $-x$ and once with $z$ replaced by $-z$, and
then subtracting.
\begin{lem} \label{Hick2}
$$
P(-x,q)P(z,q)(q)_{\infty}^2 - P(x,q)P(-z,q)(q)_{\infty}^2 =
2xP(z/x,q^2)P(xzq,q^2)(q^2;q^2)_{\infty}^2.
$$
\end{lem}
The second corollary follows just as the first, except we add
instead of subtract in the final step.
\begin{lem} \label{Hick2.5}
$$
P(-x,q)P(z,q)(q)_{\infty}^2 + P(x,q)P(-z,q)(q)_{\infty}^2 =
2P(xz,q^2)P(qz/x,q^2)(q^2;q^2)_{\infty}^2.
$$
\end{lem}
For the third corollary we subtract Lemma \ref{Hick1} with $z$
replaced by $-z$ from three times Lemma \ref{Hick1} with $x$
replaced by $-x$.
\begin{lem} \label{Hick3}
$$
\begin{aligned}
3P(-x,q)P(z,q)(q)_{\infty}^2 & - P(x,q)P(-z,q)(q)_{\infty}^2  \\
&= 2P(xz,q^2)P(zq/x,q^2)(q^2;q^2)_{\infty}^2 +
4xP(xzq,q^2)P(z/x,q^2)(q^2;q^2)_{\infty}^2.
\end{aligned}
$$
\end{lem}

Finally we record the addition theorem as stated in \cite[Eq.
(3.7)]{asd}.
\begin{lem} \label{addition}
$$
P^2(z,q)P(\zeta t,q)P(\zeta/t,q) - P^2(\zeta,q)P(zt,q)P(z/t,q) +
\zeta/t P^2(t,q)P(z\zeta,q)P(z/\zeta,q) = 0.
$$
\end{lem}

\section{Two lemmas}
The proofs of Theorems \ref{main3} and \ref{main} will follow from
identities which relate the sums $\displaystyle \Sigma(a,b)$ to
the products $P(a)$ and $P(0)$. The key steps are Lemmas
\ref{jack} and \ref{jack2} below, which are similar to Lemmas 7
and 8 in \cite{asd}.

\begin{lem} \label{jack} We have

\begin{equation} \label{jackeq}
\begin{aligned}
\sum_{n=-\infty}^{\infty} (-1)^n q^{n^2 + n} \Bigl[ \frac{\zeta^{-2n}}{1-z{\zeta^{-1}}q^{n}} +
\frac{\zeta^{2n+2}}{1-z{\zeta}q^n} \Bigr] & = \frac{\zeta (\zeta^2, q\zeta^{-2}, -1, -q)_{\infty}}{(\zeta, q\zeta^{-1}, -\zeta, -q\zeta^{-1})_{\infty}}
\sum_{n=-\infty}^{\infty} (-1)^n \frac{q^{n^2 + n}}{1-zq^n} \\
& + \frac{(\zeta, q\zeta^{-1}, \zeta^2, q\zeta^{-2}, -z, -qz^{-1}, q, q)_{\infty}}
{(z, qz^{-1}, z\zeta, qz^{-1}\zeta^{-1}, z\zeta^{-1}, q\zeta z^{-1}, -\zeta, -q\zeta^{-1})_{\infty}}.
\end{aligned}
\end{equation}

\end{lem}

\begin{proof}
This may be deduced from \cite[Eq. (4), p.236]{jackson} by making
the substitutions $a_3 = \sqrt{a_1q}$, $a_4 = z\sqrt{a_1/q}$,
$a_{10} = \sqrt{a_1q}/z$, $a_9 = -\sqrt{a_1q}$, letting $a_7$ and
$a_8$ tend to infinity, writing $\zeta = \sqrt{q/a_1}$ and
simplifying. One might also argue as in \cite[pp. 94-96]{asd} . But
the simplest way to establish the truth of \eqref{jackeq}, kindly
pointed out to us by S.H. Chan, is to observe that it is
essentially the case $r=1$, $s=3$, $a_1=-z$, $b_1=z/\zeta$, $b_2 =
z\zeta$, and $b_3=z$ of \cite[Theorem 2.1]{Ch1},
\begin{equation}
\begin{aligned}
& \frac{P(a_1,q)\cdots P(a_r,q) (q)_{\infty}^2}{P(b_1,q) \cdots
P(b_s,q)} \\ & = \frac{P(a_1/b_1,q)\cdots
P(a_r/b_1,q)}{P(b_2/b_1,q)\cdots P(b_s/b_1,q)} \sum_{n \in
\mathbb{Z}}\frac{(-1)^{(s-r)n}q^{(s-r)n(n+1)/2}}{1-b_1q^n}
\left(\frac{a_1\cdots a_r b_1^{s-r-1}}{b_2\cdots b_s}\right)^n \\
&+ \text{idem}(b_1;b_2,\dots,b_s).
\end{aligned}
\end{equation}
Here we use the notation
$$
\begin{aligned}
F(b_1,b_2,\dots ,b_m) &+  \text{idem}(b_1;b_2,\dots ,b_m) \\ &:=
F(b_1,b_2, \dots ,b_m) + F(b_2,b_1,b_3, \dots ,b_m) + \cdots +
F(b_m,b_2,\dots,b_{m-1},b_1).
\end{aligned}
$$
\end{proof}

We shall use the following specialization of Lemma \ref{jack},
which is the case $\zeta=y^a$, $z=y^b$, and $q=y^{\ell}$:

\begin{equation} \label{lem1}
\begin{aligned}
y^{2a} \Sigma(a+b, 2a) + \Sigma(b-a, -2a) -&  y^a \frac{P(2a)P(-1, y^{\ell})}{P(a)P(-y^a, y^{\ell})} \Sigma(b,0) \\
& - \frac{P(a)P(2a)P(-y^b,
y^{\ell})P(0)^2}{P(b+a)P(b-a)P(b)P(-y^a, y^{\ell})} = 0.
\end{aligned}
\end{equation}

\noindent We now define

$$
\begin{aligned}
g(z,q) & :=z\frac{P(z^2, q)P(-1,q)}{P(z,q)P(-z,q)}\Sigma(z,1,q) - z^2 \Sigma(z^2, z^2, q) \\
& -  \sideset{}{'} \sum_{n=-\infty}^{\infty} \frac{(-1)^n z^{-2n} q^{n(n+1)}}{1-q^n}
\end{aligned}
$$

\noindent and

\begin{equation} \label{g}
g(a):=g(y^a, y^{\ell})=y^a\frac{P(2a)P(-1, y^{\ell})}{P(a)P(-y^a, y^{\ell})} \Sigma(a,0) - y^{2a} \Sigma(2a, 2a) - \Sigma(0, -2a).
\end{equation}

\noindent The second key lemma is the following.

\begin{lem} \label{jack2} We have

\begin{equation} \label{part1}
2g(z,q) - g(z^2, q) + \frac{1}{2} = \frac{(q)_{\infty}^{2} P(-z^4,
q)}{P(z^4, q) P(-1, q)} + z \frac{P(-1,q)^2 (q)_{\infty}^2 P(z^2,
q)}{P(z,q)^2 P(-z, q)^2}
\end{equation}

\noindent and

\begin{equation} \label{part2}
g(z,q) + g(z^{-1}q, q)=1.
\end{equation}

\end{lem}

\begin{proof}

We first require a short computation involving $\Sigma(z, \zeta,
q)$. Note that

\begin{equation} \label{sigma}
\begin{aligned}
z^2 \Sigma(z,\zeta, q) + \zeta \Sigma(zq, \zeta, q) & =
\sum_{n=-\infty}^{\infty} (-1)^n \frac{z^2 \zeta^{n} q^{n(n+1)}}{1 - zq^n} +
\sum_{n=-\infty}^{\infty} (-1)^n \frac{\zeta^{n+1} q^{n(n+1)}}{1-zq^{n+1}} \\
& = \sum_{n=-\infty}^{\infty} \zeta^n q^{n(n-1)} \Bigl( \frac{z^2 q^{2n} -1}{1-zq^n} \Bigr) \\
& = -\sum_{n=-\infty}^{\infty} (-1)^n \zeta^n q^{n(n-1)} (1 +
zq^n)
\end{aligned}
\end{equation}

\noindent upon writing $n-1$ for $n$ in the second sum of the
first equation. Taking $\zeta=1$ yields

\begin{equation} \label{step}
\begin{aligned}
z^2 \Sigma(z,1,q) + \Sigma(zq, 1, q) & = -\sum_{n=-\infty}^{\infty} (-1)^n q^{n(n-1)} (1+ zq^n) \\
& = -z\sum_{n=-\infty}^{\infty} (-1)^n q^{n^2} \\
& = -z \prod_{r=1}^{\infty} \frac{1-q^r}{1+q^r} \\
\end{aligned}
\end{equation}

\noindent where (\ref{jtp}) was used for the last
step. Now write $g(z,q)$ in the form

$$
g(z,q)=f_{1}(z) - f_{2}(z) - f_{3}(z)
$$

\noindent where

$$
f_{1}(z):= z\frac{P(z^2, q)P(-1,q)}{P(z,q)P(-z,q)}\Sigma(z,1,q),
$$

$$
f_{2}(z):= z^2 \Sigma(z^2, z^2, q),
$$

\noindent and

$$
f_{3}(z):= \sideset{}{'} \sum_{n=-\infty}^{\infty} \frac{(-1)^n
z^{-2n} q^{n(n+1)}}{1-q^n}.
$$

\noindent By (\ref{p1}), (\ref{p2}), and (\ref{step}),

\begin{equation} \label{f1}
\begin{aligned}
f_{1}(zq) - f_{1}(z) & = \frac{P(z^2, q) P(-1, q)}{P(z,q) P(-z,q)} \prod_{r=1}^{\infty} \frac{1-q^r}{1+q^r} \\
& = 2 \sum_{n=-\infty}^{\infty} (-1)^n z^{2n} q^{n^2}. \\
\end{aligned}
\end{equation}

\noindent A similar argument as in (\ref{sigma}) yields

\begin{equation} \label{f2}
f_{2}(zq) - f_{2}(z) = \sum_{n=-\infty}^{\infty} (-1)^n z^{2n-2}
q^{n(n-1)} + \sum_{n=-\infty}^{\infty} (-1)^n z^{2n} q^{n^2}
\end{equation}

\noindent and

\begin{equation} \label{f3}
f_{3}(zq) - f_{3}(z) = -2 + \sum_{n=-\infty}^{\infty} (-1)^n
z^{-2n} q^{n(n-1)} + \sum_{n=-\infty}^{\infty} (-1)^n z^{-2n}
q^{n^2}.
\end{equation}

\noindent Adding (\ref{f2}) and (\ref{f3}), then subtracting from (\ref{f1}) gives

\begin{equation} \label{constant}
g(z,q)-g(zq, q) =-2.
\end{equation}

\noindent If we now define

$$
f(z):= 2g(z,q) - g(z^2, q) + \frac{1}{2} - \frac{(q)_{\infty}^{2}
P(-z^4, q)}{P(z^4, q) P(-1, q)} - z \frac{P(-1,q)^2 (q)_{\infty}^2
P(z^2, q)}{P(z,q)^2 P(-z, q)^2},
$$

\noindent then from (\ref{p1}), (\ref{p2}), and (\ref{constant}),
one can verify that

\begin{equation} \label{functional}
f(zq) - f(z) =0.
\end{equation}

\noindent Now, the only possible poles of $f(z)$ are simple ones
at points equivalent (under $z \to zq$) to those given by $z^4 =
1,q,q^2$, or $q^3$. For each such point $w$ it is easy to
calculate $\lim_{z\to w} (z-w)f(z)$ and see that there are, in
fact, no poles.  Hence $f(z)$ is analytic except at $z=0$, and
applying \eqref{functional} to its Laurent expansion around this
point shows that $f(z)$ is constant. Next, let us show that $f(-1)
= 0$. Since the non-constant terms in $f(z)$ have simple poles at
$z=-1$, we must consider $\lim_{z \to -1} \frac{d}{dz} (z+1)f(z)$.
We omit the computation, but mention that the term $g(z,q)$ gives
$-1/4$, $g(z^2,q)$ gives $7/8$, and the last two terms give $-1/8$
and $1$, respectively. Then $2(-1/4) - 7/8 + 1/2 + 1/8 - 1 = 0$,
and we conclude that $f(z)$ is identically $0$. This proves (\ref{part1}).

To prove (\ref{part2}), it suffices to
show, after (\ref{constant}),

\begin{equation} \label{gees}
g(z^{-1}, q) + g(z,q)=-1.
\end{equation}

\noindent Note that

\begin{equation} \label{short}
\begin{aligned}
\Sigma(z,1, q) + z^{-2} \Sigma(z^{-1}, 1, q) & =
\sum_{n=-\infty}^{\infty} (-1)^n \frac{q^{n(n+1)}}{1 - zq^n} -
z^{-1} \sum_{n=-\infty}^{\infty} (-1)^n \frac{q^{n^2}}{1-zq^{n}} \\
& = -z^{-1} \sum_{n=-\infty}^{\infty} (-1)^n q^{n^2}
\end{aligned}
\end{equation}

\noindent where we have written $-n$ for $n$ in the second sum in the first equation. Thus, by
(\ref{p1}), (\ref{p2}), and (\ref{short}), we have

\begin{equation} \label{f1z}
f_{1}(z) + f_{1}(z^{-1}) = -z^{-1} \sum_{n=-\infty}^{\infty} (-1)^n q^{n^2} \frac{P(z^2,
q)P(-1,q)}{P(z,q)P(-z,q)}.
\end{equation}

\noindent Again, a similar argument as in (\ref{short}) gives

\begin{equation} \label{f2z}
f_{2}(z) + f_{2}(z^{-1}) = -\sum_{n=-\infty}^{\infty} (-1)^n z^{2n} q^{n^2}
\end{equation}

\noindent and

\begin{equation} \label{f3z}
f_{3}(z) + f_{3}(z^{-1}) = 1 - \sum_{n=-\infty}^{\infty} (-1)^n z^{2n} q^{n^2}.
\end{equation}

\noindent Adding (\ref{f2z}) and (\ref{f3z}), then subtracting from (\ref{f1z})
yields (\ref{gees}).

\end{proof}

Letting $z=y^a$ and $q=y^{\ell}$ in Lemma \ref{jack2}, we get

\begin{equation} \label{g1}
2g(a) - g(2a) + \frac{1}{2} = \frac{P(-y^{4a}, y^{\ell})
P(0)^2}{P(4a)P(-1, y^{\ell})} + y^a \frac{P(-1, y^{\ell})^2 P(0)^2
P(2a)}{P(a)^2 P(-y^a, y^{\ell})^2}
\end{equation}

\noindent and

\begin{equation} \label{g2}
g(a) + g(q-a)=1.
\end{equation}

\noindent These two identities will be of key importance in the
proofs of Theorems \ref{main3} and \ref{main}.

\section{Proofs of Theorems \ref{main3} and \ref{main}}

We now compute the sums $\overline{S}(\ell-2m)$. The reason for
this choice is two-fold. First, we would like to obtain as simple
an expression as possible in the final formulation (\ref{final}).
Secondly, to prove Theorem \ref{main}, we will need to compute
$\overline{S}(1)$ and $\overline{S}(3)$. For $\ell=5$, we can then
choose $m=2$ and $m=1$ respectively. As this point, we follow the
idea of Section 6 in \cite{asd}. Namely, we write

\begin{equation} \label{n}
n=\ell r + m + b,
\end{equation}

\noindent where $-\infty < r < \infty$. The idea is to simplify
the exponent of $q$ in $\overline{S}(\ell-2m)$. Thus

$$
\ell n-2mn + n^2 = \ell^{2} r(r+1) + 2b\ell r + (b+m)(b-m + \ell).
$$

\noindent We now substitute (\ref{n}) into (\ref{s}) and let $b$
take the values $0$, $\pm a$, and $\pm m$. Here $a$ runs through
$1$, $2$, \dotso, $\frac{\ell-1}{2}$ where the value $a \equiv \pm
m \bmod \ell$ is omitted. As in \cite{asd}, we use the notation
$\displaystyle \sideset{}{''} \sum_{a}$ to denote the sum over
these values of $a$. We thus obtain

$$
\begin{aligned}
\overline{S}(\ell-2m) & = \sideset{}{'} \sum_{n=-\infty}^{\infty} (-1)^n \frac{q^{(\ell-2m)n + n^2}}{1-y^n} \\
& = \sum_{b} {\sideset{}{'} \sum_{r=-\infty}^{\infty}}
(-1)^{r+m+b} q^{(b+m)(b-m+\ell)} \frac{y^{\ell r(r+1) + 2br}}{1-
y^{\ell r+m+b}},
\end{aligned}
$$

\noindent where $b$ takes values $0$, $\pm a$, and $\pm m$ and the term corresponding
to $r=0$ and $b=-m$ is omitted. Thus

\begin{equation} \label{s(b)}
\begin{aligned}
\overline{S}(\ell-2m) & = (-1)^m q^{m(\ell-m)} \Sigma(m,0) + \Sigma(0, -2m) + y^{2m} \Sigma(2m, 2m) \\
& + \sideset{}{''} \sum_{a} (-1)^{m+a} q^{(a+m)(a-m + \ell)} \Bigl
\{ \Sigma(m+a, 2a) + y^{-2a} \Sigma(m-a, -2a) \Bigr \}.
\end{aligned}
\end{equation}

\noindent Here the first three terms arise from taking $b=0$, $-m$, and $m$ respectively.
We now can use (\ref{lem1}) to simplify this expression. By taking $b=m$ and dividing by
$y^{2a}$ in (\ref{lem1}), the sum of the two terms inside the curly brackets becomes

$$
y^{-a} \frac{P(2a)P(-1, y^{\ell})}{P(a)P(-y^a, y^{\ell})}
\Sigma(m,0) + y^{-2a} \frac{P(a)P(2a)P(-y^m, y^{\ell})
P(0)^2}{P(m)P(m+a)P(m-a)P(-y^a, y^{\ell})}.
$$

\noindent Similarly, upon taking $a=m$ in (\ref{g}), then the sum of the second and third terms in (\ref{s(b)}) is

$$
y^m \frac{P(2m)P(-1, y^{\ell})}{P(m)P(-y^m, y^{\ell})} \Sigma(m,0)
- g(m).
$$

In total, we have

\begin{equation} \label{final}
\begin{aligned}
\overline{S}(\ell-2m) & = -g(m) \\
& + \sideset{}{''} \sum_{a} (-1)^{m+a} q^{(a+m)(a-m+\ell)} y^{-2a}
\frac{P(a)P(2a)P(-y^m, y^{\ell})P(0)^2}{P(m)P(m+a)P(m-a)P(-y^a, y^{\ell})} \\
& + \Sigma(m,0) \Biggl\{ (-1)^m q^{m(\ell-m)} + y^m \frac{P(2m)P(-1, y^{\ell})}{P(m)P(-y^m, y^{\ell})} \\
& + \sideset{}{''} \sum_{a} (-1)^{m+a} q^{(a+m)(a-m+\ell)} y^{-a}
\frac{P(2a)P(-1, y^{\ell})}{P(a)P(-y^a, y^{\ell})} \Biggr\}.
\end{aligned}
\end{equation}

\noindent We can simplify some of the terms appearing in
(\ref{final}) as we are interested in certain values of $\ell$,
$m$, and $a$. To this end, we prove the following result. Let
$\{ \quad \}$ denote the coefficient of $\Sigma(m,0)$ in (\ref{final}).

\begin{prop} \label{brackets} If $\ell=3$ and $m=1$, then

$$
\{ \quad \} = -q^2 \frac{(q)_{\infty} (-q^{9};
q^{9})_{\infty}}{(-q)_{\infty} (q^{9}; q^{9})_{\infty}}.
$$

\noindent If $\ell=5$, $m=2$, and $a=1$, then

$$
\{ \quad \}=q^6 \frac{(q)_{\infty} (-q^{25}; q^{25})_{\infty}}{(-q)_{\infty} (q^{25}; q^{25})_{\infty}}.$$

\noindent If $\ell=5$, $m=1$, $a=2$, then

$$
\{ \quad \}=-q^4  \frac{(q)_{\infty} (-q^{25};
q^{25})_{\infty}}{(-q)_{\infty} (q^{25}; q^{25})_{\infty}}.
$$

\end{prop}

\begin{proof}
This is a straightforward application of Lemma \ref{lem6}.
\end{proof}

We are now in a position to prove Theorems \ref{main3} and
\ref{main}.  We begin with Theorem \ref{main3}.

\begin{proof}
By (\ref{gen1}), (\ref{s}), and (\ref{rels}), we have

\begin{equation} \label{gen3too}
\sum_{n=0}^{\infty} \Bigl\{ \overline{N}(0,3,n) -
\overline{N}(1,3,n) \Bigr\} q^{n}
\frac{(q)_{\infty}}{2(-q)_{\infty}}=3\overline{S}(1) +
\overline{S}(3).
\end{equation}

\noindent  By (\ref{p3}), (\ref{p4}), (\ref{final}), and
Proposition \ref{brackets} we have

\begin{equation} \label{s1too}
\overline{S}(1)= -g(1) + -q^2\Sigma(1,0) \frac{(q)_{\infty}
(-q^{9}; q^{9})_{\infty}} {(-q)_{\infty} (q^{9}; q^{9})_{\infty}}.
\end{equation}

By Lemma \ref{Sofq} we have
\begin{equation}
\overline{S}(3) = \frac{-(q)_{\infty}}{2(-q)_{\infty}} + \frac{1}{2}.
\end{equation}

We need to prove that

$$
 -3g(1) + -3q^2\Sigma(1,0) \frac{(q)_{\infty}
(-q^{9}; q^{9})_{\infty}} {(-q)_{\infty} (q^{9}; q^{9})_{\infty}}
- \frac{(q)_{\infty}}{2(-q)_{\infty}} + \frac{1}{2}  = \Bigl \{
r_{01}(0) q^0 + r_{01}(1)q + r_{01}(2)q^2  \Bigr \}
\frac{(q)_{\infty}}{2(-q)_{\infty}}.
$$

\noindent We now multiply the right hand side of the above
expression using Lemma \ref{lem6} and the $r_{01}(d)$ from Theorem
\ref{main} (recall that $r_{01}(d)$ is just $R_{01}(d)$ with $q$
replaced by $q^{3}$).  We then equate coefficients of powers of
$q$ and verify the resulting identities.  The only power of $q$
for which the resulting equation does not follow easily upon
cancelling factors in infinite products is the constant term. We
obtain
$$
-3g(1) + \frac{1}{2} =
\frac{(q^9;q^9)_{\infty}^3(-q^3;q^3)_{\infty}}{2(q^3;q^3)_{\infty}(-q^9;q^9)_{\infty}^3}
- 4q^3
\frac{(-q^9;q^9)_{\infty}^3(q^{18};q^{18})_{\infty}^3}{(q^6;q^6)_{\infty}(-q^3;q^3)_{\infty}}.
$$
But this follows from \eqref{g1} and some simplification.  This
then completes the proof of Theorem \ref{main3}.
\end{proof}

We now turn to Theorem \ref{main}.

\begin{proof}
We begin with the rank differences $R_{12}(d)$.  By (\ref{gen1}),
(\ref{s}), and (\ref{rels}), we have

\begin{equation} \label{gen3}
\sum_{n=0}^{\infty} \Bigl\{ \overline{N}(1,5,n) -
\overline{N}(2,5,n) \Bigr\} q^{n}
\frac{(q)_{\infty}}{2(-q)_{\infty}}=-\overline{S}(1) -
3\overline{S}(3)
\end{equation}

\noindent and by (\ref{p3}), (\ref{p4}), (\ref{final}), and Proposition \ref{brackets},

\begin{equation} \label{s1}
\overline{S}(1)= -g(2) + qy\Sigma(2,0) \frac{(q)_{\infty}
(-q^{25}; q^{25})_{\infty}} {(-q)_{\infty} (q^{25};
q^{25})_{\infty}} - q^{2} \frac{(q^{25}; q^{25})_{\infty}^2
(-q^{10}, -q^{15}; q^{25})_{\infty}} {(q^{10}, q^{15};
q^{25})_{\infty} (-q^{5}, -q^{20}; q^{25})_{\infty}}
\end{equation}

\noindent and

\begin{equation} \label{s3}
\overline{S}(3) = -g(1) - q^4 \Sigma(1,0) \frac{(q)_{\infty}
(-q^{25}; q^{25})_{\infty}} {(-q)_{\infty} (q^{25};
q^{25})_{\infty}} + q^3 \frac{(q^{25}; q^{25})_{\infty}^2 (-q^{5},
-q^{20}; q^{25})_{\infty}} {(q^{5}, q^{20}; q^{25})_{\infty}
(-q^{10}, -q^{15}; q^{25})_{\infty}}.
\end{equation}

\noindent By (\ref{gen3}), (\ref{s1}), and (\ref{s3}), we need to prove

$$
\begin{aligned}
& g(2) - qy\Sigma(2,0) \frac{(q)_{\infty} (-q^{25};
q^{25})_{\infty}} {(-q)_{\infty} (q^{25}; q^{25})_{\infty}} +
q^{2} \frac{(q^{25}; q^{25})_{\infty}^2 (-q^{10}, -q^{15};
q^{25})_{\infty}}
{(q^{10}, q^{15}; q^{25})_{\infty} (-q^{5}, -q^{20}; q^{25})_{\infty}} \\
& +3g(1) + 3q^4 \Sigma(1,0) \frac{(q)_{\infty} (-q^{25};
q^{25})_{\infty}} {(-q)_{\infty} (q^{25}; q^{25})_{\infty}} -3q^3
\frac{(q^{25}; q^{25})_{\infty}^2 (-q^{5}, -q^{20};
q^{25})_{\infty}}
{(q^{5}, q^{20}; q^{25})_{\infty} (-q^{10}, -q^{15}; q^{25})_{\infty}} \\
& = \Bigl \{ r_{12}(0) q^0 + r_{12}(1)q + r_{12}(2)q^2 +
r_{12}(3)q^3 + r_{12}(4)q^4 \Bigr \}
\frac{(q)_{\infty}}{2(-q)_{\infty}}.
\end{aligned}
$$

\noindent We now multiply the right hand side of the above
expression using Lemma \ref{lem6} and the $R_{12}(d)$ from Theorem
\ref{main}, equating coefficients of powers of $q$.  The
coefficients of $q^{0}$, $q^{1}$, $q^{2}$, $q^{3}$, $q^{4}$ give
us, respectively,

\begin{equation} \label{check0}
g(2) + 3g(1) = y \frac{(q^{25}; q^{25})_{\infty}^2}{(q^{15},
q^{20}, q^{30}, q^{35}; q^{50})_{\infty}} + 4y \frac{(q^{10},
q^{15}, q^{35}, q^{40}; q^{50})_{\infty} (q^{50};
q^{50})_{\infty}^2} {(q^{20}, q^{30}; q^{50})_{\infty}^2 (q^5,
q^{45}; q^{50})_{\infty}},
\end{equation}

\begin{equation} \label{check1}
y \frac{(q^{50}; q^{50})_{\infty} (q^{15}, q^{35}, q^{50};
q^{50})_{\infty}} {(q^{15}, q^{20}, q^{30}, q^{35};
q^{50})_{\infty}} = y \frac{(q^{50}; q^{50})_{\infty} (q^{5},
q^{45}, q^{50}; q^{50})_{\infty}}{(q^5, q^{20}; q^{25})_{\infty}},
\end{equation}

\begin{equation} \label{check2}
\frac{(q^{25}; q^{25})_{\infty}^2 (-q^{10}, -q^{15};
q^{25})_{\infty}} {(q^{10}, q^{15}; q^{25})_{\infty} (-q^5,
-q^{20}; q^{25})_{\infty}} = \frac{(q^{25};
q^{25})_{\infty}^2}{(q^5, q^{20}; q^{25})_{\infty}} - 2y
\frac{(q^{50}; q^{50})_{\infty}^2 (q^{5}, q^{45};
q^{50})_{\infty}}{(q^{10}, q^{15}; q^{25})_{\infty}},
\end{equation}

\begin{equation} \label{check3}
\begin{aligned}
\frac{3(q^{25}; q^{25})_{\infty}^2 (-q^{5}, -q^{20};
q^{25})_{\infty}} {(q^5, q^{20}; q^{25})_{\infty} (-q^{10},
-q^{15}; q^{25})_{\infty}} & = \frac{(q^{25};
q^{25})_{\infty}^2}{(q^{10}, q^{15}; q^{25})_{\infty}} +
\frac{2(q^{50}; q^{50})_{\infty}^2 (q^{15}, q^{35}; q^{50})}{(q^{5}, q^{20}; q^{25})_{\infty}} \\
& + 4y\frac{(q^{10}, q^{40}; q^{50})_{\infty} (q^{50};
q^{50})_{\infty}^2}{(q^{20}, q^{30}; q^{50})_{\infty}^2},
\end{aligned}
\end{equation}

\begin{equation} \label{check4}
\begin{aligned}
\frac{(q^{10}, q^{40}, q^{50}; q^{50})_{\infty} (q^{25};
q^{25})_{\infty}} {(q^{20}, q^{30}; q^{50})_{\infty}^2 (q^5;
q^{45}; q^{50})_{\infty} (-q^{25}; q^{25})_{\infty}} & =
\frac{(q^{50}; q^{50})_{\infty}^2 (q^{15}, q^{35}; q^{50})_{\infty}}{(q^{10}, q^{15}; q^{25})_{\infty}} \\
& + y\frac{(q^{50}; q^{50})_{\infty}^2 (q^5, q^{45};
q^{50})_{\infty}}{(q^{15}, q^{20}, q^{30}, q^{35};
q^{50})_{\infty}}.
\end{aligned}
\end{equation}

Equation \eqref{check1} is immediate.  Upon clearing denominators
in \eqref{check2}-\eqref{check4} and simplifying, we see
that \eqref{check2} is equivalent to the case $(x,z,q) =
(-q^5,-q^{10},q^{25})$ of Lemma \ref{Hick2}, \eqref{check3} is the
case $(x,z,q) = (q^5,q^{10},q^{25})$ of Lemma \ref{Hick3}, and
\eqref{check4} follows from the case $(z,\zeta,t,q) =
(q^{20},q^{10},q^{5},q^{50})$ of Lemma \ref{addition}.

As for \eqref{check0}, let us take $a=1$ and $a=2$ in \eqref{g1},
and then replace $g(4)$ by $1-g(1)$ using \eqref{g2}.  This gives
$$
\begin{aligned}
3g(1) + g(2) &=
\frac{(-q^5,-q^{20},q^{25},q^{25};q^{25})_{\infty}}{2(q^5,q^{20},-q^{25},-q^{25};q^{25})_{\infty}}
-
\frac{4y(q^{15},q^{35},q^{50},q^{50};q^{50})_{\infty}}{(q^{10},q^{25},q^{25},q^{40};q^{50})_{\infty}}
\\
&-
\frac{(-q^{10},-q^{15},q^{25},q^{25};q^{25})_{\infty}}{2(q^{10},q^{15},-q^{25},-q^{25};q^{25})_{\infty}}
+
\frac{4y^2(q^{5},q^{45},q^{50},q^{50};q^{50})_{\infty}}{(q^{20},q^{25},q^{25},q^{30};q^{50})_{\infty}}.
\end{aligned}
$$
Now after making a common denominator in the first and third
terms, we may apply the case $(x,z,q) = (q^5,q^{10},q^{25})$ of
Lemma \ref{Hick2} to these two terms, the result being precisely
the first term in \eqref{check0}.  For the second and fourth
terms, we make a common denominator and multiply top and bottom by
$(q^5,q^{20};q^{25})_{\infty}$.  Then the case $(z, \zeta,t,q) =
(q^{20},q^{10},q^{5},q^{50})$ of Lemma \ref{addition} applies and
we obtain the second term in \eqref{check0}.

We now turn to the rank differences $R_{02}(d)$, proceeding as
above. Again by (\ref{gen1}), (\ref{s}), and (\ref{rels}), we have

\begin{equation} \label{gen4}
\sum_{n=0}^{\infty} \Bigl\{ \overline{N}(0,5,n) -
\overline{N}(2,5,n) \Bigr\} q^{n}
\frac{(q)_{\infty}}{2(-q)_{\infty}}=-\overline{S}(5) +
2\overline{S}(1) + \overline{S}(3).
\end{equation}

\noindent By Lemma \ref{Sofq} (with $\ell=5$), (\ref{gen4}),
(\ref{s1}), and (\ref{s3}), it suffices to prove

$$
\begin{aligned}
&  \frac{-(q)_{\infty}}{2(-q)_{\infty}} + \frac{1}{2} - 2g(2) +
2qy\Sigma(2,0) \frac{(q)_{\infty} (-q^{25}; q^{25})_{\infty}}
{(-q)_{\infty} (q^{25}; q^{25})_{\infty}} - 2q^{2} \frac{(q^{25};
q^{25})_{\infty}^2 (-q^{10}, -q^{15}; q^{25})_{\infty}}
{(q^{10}, q^{15}; q^{25})_{\infty} (-q^{5}, -q^{20}; q^{25})_{\infty}} \\
& -g(1) - q^4 \Sigma(1,0) \frac{(q)_{\infty} (-q^{25};
q^{25})_{\infty}} {(-q)_{\infty} (q^{25}; q^{25})_{\infty}} + q^3
\frac{(q^{25}; q^{25})_{\infty}^2 (-q^{5}, -q^{20};
q^{25})_{\infty}}
{(q^{5}, q^{20}; q^{25})_{\infty} (-q^{10}, -q^{15}; q^{25})_{\infty}} \\
& = \Bigl \{ r_{02}(0) q^0 + r_{02}(1)q + r_{02}(2)q^2 +
r_{02}(3)q^3 + r_{02}(4)q^4 \Bigr \}
\frac{(q)_{\infty}}{2(-q)_{\infty}}.
\end{aligned}
$$

\noindent Again, equating coefficients of powers of $q$ yields the
following identities.

\begin{equation} \label{check5}
\begin{aligned}
\frac{1}{2} -2g(2) - g(1) & = \frac{1}{2} \frac{(-q^{10}, -q^{15};
q^{25})_{\infty} (q^{25}; q^{25})_{\infty}^2} {(q^{10}, q^{15};
q^{25})_{\infty} (-q^{25}; q^{25})_{\infty}^2} - 2y \frac{(q^{10},
q^{40}; q^{50})_{\infty} (q^{15}, q^{35}; q^{50})_{\infty}
(q^{50}; q^{50})_{\infty}^2}
{(q^{20}, q^{30}; q^{50})_{\infty}^2 (q^5, q^{45}; q^{50})_{\infty}} \\
& + 2y\frac{(q^{20}, q^{30}; q^{50})_{\infty} (q^5, q^{45};
q^{50})_{\infty} (q^{50}; q^{50})_{\infty}^2} {(q^{10}, q^{40};
q^{50})_{\infty}^2(q^{15},q^{35};q^{50})_{\infty}},
\end{aligned}
\end{equation}

\begin{equation} \label{check6}
\frac{(q^{20}, q^{30}, q^{50}; q^{50})_{\infty} (q^{25};
q^{25})_{\infty}}
{(q^{10},q^{40}; q^{50})_{\infty}^2 (q^{15},q^{35};q^{50})_{\infty} (-q^{25}; q^{25})_{\infty}} \\
= \frac{(-q^{10}, -q^{15}; q^{25})_{\infty} (q^{25};
q^{25})_{\infty} (q^{15}, q^{35}, q^{50}; q^{50})_{\infty}}
{(q^{10}, q^{15}; q^{25})_{\infty} (-q^{25}; q^{25})_{\infty}},
\end{equation}

\begin{equation} \label{check7}
\begin{aligned}
& \frac{(q^{25}; q^{25})_{\infty}^2 (-q^{10}, -q^{15};
q^{25})_{\infty}}
{(-q^5, -q^{20}; q^{25})_{\infty} (q^{10}, q^{15}; q^{25})_{\infty}} \\
& = \frac{(q^{50}; q^{50})_{\infty}^2  (q^{20}, q^{30};
q^{50})_{\infty}} {(q^{10}, q^{40}; q^{50})_{\infty}^2}
- y\frac{(q^{50}; q^{50})_{\infty}^2 (q^{5}, q^{45}; q^{50})}{(q^{10}, q^{15}; q^{25})_{\infty}}, \\
\end{aligned}
\end{equation}

\begin{equation} \label{check8}
\frac{(q^{25}; q^{25})_{\infty}^2 (-q^{5}, -q^{20};
q^{25})_{\infty}} {(-q^{10}, -q^{15}; q^{25})_{\infty} (q^5,
q^{20}; q^{25})_{\infty}} = \frac{(q^{25};
q^{25})_{\infty}^2}{(q^{10}, q^{15}; q^{25})_{\infty}} + 2y
\frac{(q^{50}; q^{50})_{\infty}^2 (q^{10}, q^{40};
q^{50})_{\infty}}{(q^{20}, q^{30}; q^{50})_{\infty}^2},
\end{equation}

\begin{equation} \label{check9}
\begin{aligned}
\frac{(q^{10}, q^{40}, q^{50}; q^{50})_{\infty} (q^{25};
q^{25})_{\infty}} {(q^{20}, q^{30}; q^{50})_{\infty}^2 (q^5;
q^{45}; q^{50})_{\infty} (-q^{25}; q^{25})_{\infty}} &+
\frac{(-q^{10}, -q^{15}; q^{25})_{\infty} (q^{25};
q^{25})_{\infty} (q^5,q^{45},q^{50};q^{50})_{\infty}} {(q^{10},
q^{15};
q^{25})_{\infty} (-q^{25}; q^{25})_{\infty}} \\
&= 2\frac{(q^{50}; q^{50})_{\infty}^2 (q^{15}, q^{35};
q^{50})_{\infty}}{(q^{10}, q^{15}; q^{25})_{\infty}}.
\end{aligned}
\end{equation}
Now, \eqref{check6} is immediate.  After clearing denominators and
simplifying, \eqref{check7} follows from the case $(z,\zeta,t,q) =
(q^{20},q^{15},q^{10},q^{50})$ of Lemma \ref{addition} and
\eqref{check8} is the case $(x,z,q) = (q^{5},q^{10},q^{25})$ of
Lemma \ref{Hick2}.  For \eqref{check9}, we simplify the first term
and apply the case $(x,z,q) = (q^5,q^{10},q^{25})$ of Lemma
\ref{Hick2.5}.

As for \eqref{check5}, taking the case $a=2$ of \eqref{g1}
together with an application of \eqref{g2} gives
$$
\frac{1}{2} -2g(2) - g(1)  = \frac{(-q^{10}, -q^{15};
q^{25})_{\infty} (q^{25}; q^{25})_{\infty}^2} {2(q^{10}, q^{15};
q^{25})_{\infty} (-q^{25}; q^{25})_{\infty}^2} -
\frac{4y^2(q^5,q^{45};q^{50})_{\infty}(q^{50};q^{50})_{\infty}^2}
{(q^{20},q^{30};q^{50})_{\infty}(q^{25};q^{50})_{\infty}^2}.
$$
Now the first terms of the above equation and \eqref{check5} match
up, while after some simplification of the final two terms of
\eqref{check5} we may apply the case $(x,z,q) =
(q^5,q^{10},q^{25})$ of Lemma \ref{Hick2} to obtain the final term
above.
\end{proof}

\section*{Acknowlegements}
The second author would like to thank the Institut des Hautes {\'E}tudes
Scientifiques and the Max-Planck-Institut f{\"u}r Mathematik
for their hospitality and support during the preparation
of this paper.

\end{document}